\documentclass[pdflatex]{sn-jnl}

\usepackage{amsrefs}

\usepackage{graphicx}%
\usepackage{multirow}%
\usepackage{amsmath,amssymb,amsfonts}%
\usepackage{amsthm}%
\usepackage{mathrsfs}%
\usepackage[title]{appendix}%
\usepackage{xcolor}%
\usepackage{textcomp}%
\usepackage{manyfoot}%
\usepackage{booktabs}%
\usepackage{algorithm}%
\usepackage{algorithmicx}%
\usepackage{algpseudocode}%
\usepackage{listings}%
\usepackage{xspace}
\usepackage{amsrefs}
\usepackage{bm}
\usepackage{tikz}
\usetikzlibrary{arrows,arrows.meta,matrix}
\usepackage{enumerate}
\makeatletter
\renewcommand{\proofname}{Proof.}
\newenvironment{myspiproof}[1][\proofname]{\par\removelastskip
  \pushQED{\qed}%
  \normalfont \topsep7.5\p@\@plus7.5\p@\relax%
  \trivlist%
  \item[\hskip\labelsep%
        \itshape%
    #1\@addpunct{}]\ignorespaces%
}{%
  \popQED\endtrivlist\@endpefalse%
}%

\newtheoremstyle{mythmstyleone}
{18pt plus2pt minus1pt}
{18pt plus2pt minus1pt}
{\normalsize\itshape}
{0pt}
{\bfseries}
{}
{.5em}
{\thmname{#1}\thmnumber{\@ifnotempty{#1}{ }\@upn{#2}}%
 . \thmnote{ {\the\thm@notefont(#3)}}}
\newtheoremstyle{mythmstylethree}
{18pt plus2pt minus1pt}
{18pt plus2pt minus1pt}
{\normalfont}
{0pt}
{\bfseries}
{}
{.5em}
{\thmname{#1}\thmnumber{\@ifnotempty{#1}{ }\@upn{#2}}%
 . \thmnote{ {\the\thm@notefont(#3)}}}

\makeatother

\theoremstyle{mythmstyleone}%
\newtheorem{theorem}{Theorem}[section]
\newtheorem{lemma}[theorem]{Lemma}
\newtheorem{corollary}[theorem]{Corollary}
\newtheorem{problem}[theorem]{Problem}
\newtheorem*{Claim}{Claim}
\newtheorem*{mlemma}{Main Lemma}

\theoremstyle{mythmstylethree}%
\newtheorem{definition}[theorem]{Definition}%

\def\<{\left\langle}
\def\>{\right\rangle}
\newcommand*\subs{\subset}
\newcommand*\subsq{\subseteq}

\newcommand*\nb{Noetherian base\xspace}
\newcommand*\ob{outer base\xspace}
\newcommand*\nob{Noetherian outer base\xspace}
\newcommand*\nobs{Noetherian outer bases\xspace}
\newcommand*\lno{$\lambda$-Noetherian\xspace}
\newcommand*\lnob{$\lambda$-Noetherian outer base\xspace}

\newcommand{\Nt}{\operatorname{Nt}}
\newcommand{\mc}{\mathcal}
\newcommand{\doublearrow}{\xrightarrow[]{}\mathrel{\mkern-14mu}\rightarrow}

\newcommand*\set[1]{\left\{#1\right\}}
\newcommand*\abs[1]{\left\vert #1\right\vert}

\DeclareMathOperator{\RO}{RO}

\DeclareMathOperator{\tr}{tr}
\newcommand*\wt[1]{\widetilde{#1}}
\newcommand*\sm{\setminus}
\newcommand*\0{\emptyset}

\newcommand*\cl[1]{\overline{#1}}
\DeclareMathOperator{\lu}{lu}
\DeclareMathOperator{\cf}{cf}
\DeclareMathOperator{\inte}{int}

\newcommand\ds{\displaystyle}

\newcommand{\setm}{\setminus}
\newcommand{\empt}{\emptyset}

\def\<{\left\langle}
\def\>{\right\rangle}
\def\br#1;#2;{\bigl[ {#1} \bigr]^ {#2} }

\newcommand{\alo}{{\aleph_{\omega}}}

\newcommand{\deltop}[1]{{#1}_{\delta}}

\newcommand{\dom}{\operatorname{dom}}
\newcommand{\ran}{\operatorname{ran}}

\begin{document}

\title[Noetherian Properties]{Noetherian Properties, Large Cardinals, and  
Independence  Around  $\aleph_{\omega}$ }


\author*[1]{\fnm{Lajos} \sur{Soukup}}\email{soukup@renyi.hu}

\author[2]{\fnm{Zolt\'an} \sur{Szentmikl\'ossy}}\email{szentmiklossy@renyi.hu}

\affil*[1]{\orgdiv{Department of  Set-Theory, Topology and Logic}, \orgname{HUN-REN Alfr\'ed R\'enyi Institute of Mathematics}, 
\orgaddress{\street{Re\'altanoda utca, 13--15}, 
\city{Budapest}, \postcode{H1053},  \country{Hungary}}}

\affil[2]{\orgdiv{Department of Analysis}, \orgname{Eötvös Lor{\'a}nt University}, 
\orgaddress{\street{P{\'a}zm{\'a}ny P{\'e}ter s{\'e}t{\'a}ny 1/A}, \city{Budapest}, \postcode{H1117},  \country{Hungary}}}


\abstract{
    A base of a topological space is called {\em Noetherian } 
iff it does not contain an
infinite strictly $\subseteq$-increasing chain.

We show that  minimal  
cardinality of  a regular spaces without a \nb is the 
first strongly inaccessible cardinal, answering a question from the 1980s.

We also study 
the {\em Noetherian type} of a topological space $X$, denoted by $\Nt(X)$,  defined 
as
the least cardinal $\kappa$ such that $X$ has a base $\mc B$ with 
$|\{B'\in \mc B: B\subs B'\}|<\kappa$ for each $B\in \mc B$. 
The behavior of the Noetherian type  under the  
$G_\delta$-modification was investigated by Milovich and Spadaro.
A central question, posed by them, is whether 
the Noetherian type of the $G_{\delta}$-modification of the space 
    $D(2)^{\aleph_\omega}$ is $\omega_1$. This statement,
 denoted (Nt), 
is known to be independent of ZFC + GCH: it holds under 
``GCH + $\square_{\aleph_\omega}$'', but fails under 
``GCH + $(\aleph_{\omega+1}, \aleph_\omega) \doublearrow (\aleph_1, \aleph_0)$''.

We place this phenomenon in a broader context by identifying similar independence phenomena for several topological and 
combinatorial principles. These include: (wFN) the weak Freese-Nation property of $[\aleph_{\omega}]^{\omega}$; (SAT) the existence of a saturated MAD family
in $[\aleph_{\omega}]^{\omega}$; (HnT) the existence of  an ${\omega}$-homogeneous, but not ${\omega}$-transitive permutation 
group on $\aleph_{\omega}$; and (SPL) the existence of a  
 countably compact, locally countable, and ${\omega}$-fair regular space of cardinality $\aleph_{\omega+1}$.
 Assuming GCH, we analyze the logical relationships between 
 these principles and show, for example, that
 SPL implies wFN, which in turn implies both SAT, HnT and Nt, while SAT does not imply Nt.
}



\maketitle

\section{Introduction}
The main question addressed in the first part of this paper is to identify the  topological spaces that possess a so-called {\em Noetherian base}, defined as an open base without a proper infinite $\subset$-increasing sequence.
This notion was introduced by Lindgren and Nyikos in 
\cite{LiNyi76}*{Definition 3.1}. 
Spaces with a Noetherian base have been considered by  Arhangel'skii,  Choban, Förrster,  Grabner,  Gruenhage,  Lindgren,  
Malykhin,  Nyikos and  Peregudov   
(see \cites{ArCh12,FoGr82,Gr83,LiNyi76,Ma81,Pe97}).

It is easy to see that every metric space has a Noetherian base. 
However, as noted in \cite{Ta85}, van Douwen --- in an unpublished privately circulated  note --- proved 
that an ordinal ${\alpha}$ with the order topology has a Noetherian base if and only if  ${\alpha}$ is less than the 
first  strongly inaccessible cardinal (for a proof, see \cite{Ta85}*{Theorem  2.8}). In \cite{Ta88}, the authors present an example of a $T_1$-space which does not have a Noetherian base in ZFC, and  they mentioned  as 
interesting unsolved problem  whether there are topological spaces 
with stronger separation axioms that lack a Noetherian base

We aim to demonstrate that the notion of a strongly inaccessible cardinal plays a crucial role in the general case --- 
specifically, in the  following theorem:

\begin{theorem}\label{tm:noether-base}
    If $X$ is a $T_3$ space and $\abs X$ is less than the first strongly inaccessible cardinal, then $X$ has a \nb.
   \end{theorem}
   
However, the following question remained open. 
   \begin{problem}\label{pr:t2}
   Is there a $T_2$-space which does not have a Noetherian base in ZFC?
   \end{problem}

 \medskip

The Noetherian type of a space is an order-theoretic variant on weight that was introduced by Peregudov (\cite{Pe84}), although 
investigations on it go back to Malykhin (\cite{Ma81}), and Peregudov and Shapirovskii (\cite{PeSh76}).

\begin{definition}
   The {\em Noetherian type} of a space $X$, denoted by $\Nt(X)$, is 
   the least cardinal $\kappa$ such that $X$ has a base $\mc B$ such that 
   $|\{B'\in \mc B: B\subs B'\}|<\kappa$ for each $B\in \mc B$. 
\end{definition}

The {\em $G_{\delta}$-modification} of a topological space $X$,
denoted by $X_{\delta}$, is the topology generated by all $G_{\delta}$-subsets of 
$X$.

The behavior of the Noetherian type  under the  
$G_\delta$-modification was investigated by Milovich and Spadaro.
A central question posed  by them is whether 
$\Nt((D(2)^{\aleph_\omega})_\delta)=\omega_1$. 
This statement, denoted (Nt),  is known to be independent of 
ZFC + GCH: it holds under ``GCH + $\square_{\aleph_\omega}$'' by \cite{KoMiSp14}, 
but fails under 
``GCH +  $(\aleph_{\omega+1}, \aleph_\omega) \doublearrow (\aleph_1, \aleph_0)$'' 
by \cite{soukup2010notenoetheriantypespaces}. We note that 
 this partition relation constitutes a large cardinal assumption, as shown by \cite{LeMaSh90}. 
  
 In the second part of this paper we place this phenomenon in a broader context by identifying similar independence phenomena for several topological and 
      combinatorial principles. These include: (wFN) the weak Freese-Nation property of $[\aleph_{\omega}]^{\omega}$; (SAT) the existence of a saturated MAD family
      in $[\aleph_{\omega}]^{\omega}$; (HnT) the existence of  an ${\omega}$-homogeneous, but not ${\omega}$-transitive permutation 
      group on $\aleph_{\omega}$; and (SPL) the existence of a  
       countably compact, locally countable, and ${\omega}$-fair regular space of cardinality $\aleph_{\omega+1}$.

       In  Theorem \ref{tm:main-implication}, assuming GCH, we analyze the logical relationships between 
       these principles and show, for example, that
       SPL implies wFN, which in turn implies both SAT, HnT and Nt, while SAT does not imply Nt.

\section{Noetherian bases}

We begin with some definitions and notation.

A partially ordered set $\< P, \le \>$ is called {\em Noetherian} if 
it contains no infinite strictly increasing sequence; that is, it is ``upward well-founded''.

Let $\lambda$ be a cardinal. We say that a partially ordered set $\< P, \le \>$ is {\em \lno} 
if it can be partitioned into at most $\lambda$ many Noetherian subsets. 
That is, there exists a function $f: P \to \lambda$ such that 
for each $\xi \in \lambda$, the poset
\[
\< f^{-1}(\set{\xi}), \le \restriction f^{-1}(\set{\xi}) \>
\]
is Noetherian.

\begin{lemma}\label{MonNoether}
   Let $\<P,\le\>$ and $\<Q,\le\>$ be partially ordered sets, and let ${\varphi}: P \to Q$ be a monotonic function. If $Q$ is \lno and, for every $q \in Q$, the fiber ${\varphi}^{-1}(\set{q})$ is also \lno, then $P$ is \lno as well.
   \end{lemma}
   
   \begin{proof}
   Let $g: Q \to \lambda$ and, for each $q \in Q$, let $g_q: {\varphi}^{-1}(\set{q}) \to \lambda$ be functions witnessing that both $Q$ and the fibers are \lno, respectively. Define a function $h: P \to \lambda \times \lambda$ by
   \[
   h(p) = \<g_{{\varphi}(p)}(p), g({\varphi}(p))\>.
   \]
   We claim that $h$ decomposes $P$ into Noetherian pieces.
   
 Fix  $\<{\xi},{\eta}\>\in {\lambda}\times {\lambda}$ and  consider an arbitrary infinite increasing sequence $\<p_n : n \in \omega\>$ in $P$ such that 
   $h(p_n) = \<{\xi},{\eta}\>$ for all $n \in \omega$.
   
   The sequence $\<{\varphi}(p_n) : n \in \omega\>$ is increasing in $Q$, since ${\varphi}$ is monotonic. Furthermore,  we have 
   $g({\varphi}(p_n)) = {\eta}$ for all $n\in {\omega}$ because $h(p_n) = \<{\xi},{\eta}\>$.
Since $g$ witness that $Q$ is \lno, it follows 
that    the sequence $\<{\varphi}(p_n) : n \in \omega\>$ is eventually constant. That is, 
   there exists $M<{\omega}$ and $q\in Q$
   such that   ${\varphi}(p_m) = q$ for all  $m\ge M$. 
   
Hence $\{p_m:M\le m<{\omega}\}\subs {\varphi}^{-1}\{q\}$ and 
$g_q(p_m)={\xi}$ for all $M\le m<{\omega}$ because $h(p_m)=\<{\xi},{\eta}\>$.
Since  $g_q$ witnesses that $\varphi^{-1}(\{q\})$ is \lno, the sequence $\<p_m : m < \omega\>$ must be eventually constant.
   
   This shows that $h$ indeed decomposes $P$ into Noetherian pieces.
   \end{proof}
   
\newcommand{\klfamily}[2]{P_{#1}^{#2}}

   \begin{lemma}\label{LuNoether}
For any cardinals $\kappa$ and $\lambda$, the poset 
   $\klfamily{{\kappa}}{<{\lambda}}=\<[\kappa]^{<\lambda},\subseteq\>$ is $2^{<\lambda}$-Noetherian.
   \end{lemma}
   
   \begin{proof}  It is a folklore result that 
       $\klfamily{{\kappa}}{<{\lambda}}$ is the union of $2^{<{\lambda}}$-many antichains.  
      For a short proof of that fact,  see the argument presented in the  MathOverflow post \cite{bof}.
   \end{proof}
   
   Let $\<X, \tau\>$ be a topological space, and let ${\mc  G} \subseteq \tau$. We say that {\em ${\mc  G}$ has the Noetherian property}, or briefly that {\em ${\mc  G}$ is Noetherian}, if the poset $\<{\mc  G}, \subseteq\>$ is Noetherian.
   
   \smallskip\noindent
   Similarly, we say that ${\mc  G}$ is {\em \lno} if $\<{\mc  G}, \subseteq\>$ is \lno.
   
   \smallskip\noindent
   A {\em \nb}  of a space $X$ is a base that is Noetherian.
   
   \smallskip\noindent
   Now let $\<X, \tau\>$ be a topological space, $Y \subseteq X$, and ${\mc  B} \subseteq \tau$. We say that
   {\em ${\mc  B}$ is an outer base of $Y$ in $X$} if, for every $p \in Y$, the set 
   \[
   \set{G \in {\mc  B} : p \in G}
   \]
   is a neighborhood base of $p$ in $X$.
   
   \medskip\noindent
   For a point $p \in X$, define
   \[
   \varrho(p, X) = \min\set{\rho({G}) : p \in G \in \tau},
   \]
   where $\RO(G)$ denotes the family of regular open subsets of $G$, and $\rho(G)=|\RO(G)|$.

   The following lemma is straightforward,
   \begin{lemma}\label{NoUnion}
   Let $X$ be a topological space $X$ that is the disjoint union of open sets, each of which has a Noetherian base. Then $X$ has a Noetherian base.
   \end{lemma}

\subsection*{$\varrho$-homogeneous spaces}
\begin{lemma}\label{lem1}
 Let $\<X,\tau\>$ be a topological space.
\begin{enumerate}[(a)]
 \item If $Y\subseteq X$ and $\<Y,\tau\restriction Y\>$ has a \nb, then $Y$ has a \nob in $X$.
\item Every point $p\in X$ has a Noetherian neighborhood base.
\item If both $Y$ and $Z$ have \nobs, then their union $Y\cup Z$ also has a \nob.
\end{enumerate}

\end{lemma}

\begin{proof}\ 
 
\noindent(a)
 Let $\widetilde{{\mc  B}}$ be a \nb for the subspace $Y$. For each $U\in \widetilde{{\mc  B}}$, define
\begin{displaymath}
   {\mc  B}_U = \set{G\in \tau: G\cap Y = U}, 
\end{displaymath}
fix an enumeration 
$$   \mc B_U=\set{G_\xi:\xi \in \abs{{\mc  B}_U}},$$
and let 
\[
 {\mc  B}'_U = \set{G_\xi\in{\mc  B}_U: \forall\eta<\xi \; (G_\eta \not\subset G_\xi)}.
\]
The family $\mathcal{B}'_U$ is Noetherian: any strictly increasing sequence in $\mathcal{B}'_U$ would yield a strictly decreasing sequence of indices, which is impossible. 

Moreover, for every $G' \in \tau$ with $U \subseteq G'$, there exists $G \in \mathcal{B}'_U$ such that $G \subseteq G'$.
Hence, the family
\[
 {\mc  B} = \bigcup\set{{\mc  B}'_U: U\in\widetilde{{\mc  B}} }
\]
forms an  outer-base of $Y$ in $X$. Define the map  
$f:{\mc  B}\to\widetilde{{\mc  B}}$ by 
\[
 f(G) = G\cap Y.
\]
Since the function $f$ is monotonic, each fiber $f^{-1}(\{U\}) = \mathcal{B}'_U$ is Noetherian and $\widetilde{\mathcal{B}}$ is Noetherian, it follows from Lemma~\ref{MonNoether} (with $\lambda = 1$) that $\mathcal{B}$ is also Noetherian.

\medskip \noindent (b)
 Apply the previous statement to the singleton $Y=\set{p}$.

 \medskip \noindent (c)
 Let  ${\mc  B}_Y$ and ${\mc  B}_Z$ be  \nobs of $Y$ and $Z$, respectively. Then
${\mc  B} = {\mc  B}_Y \cup {\mc  B}_Z$ is a \nob of $Y\cup Z$.
\end{proof}

\begin{corollary}\label{cor}Any discrete subset of a topological space 
   $X$ admits an outer Noetherian base.
\end{corollary}

\begin{mlemma}
 Let $X$ be a $T_3$ space, and let ${\mc  B} \subs \RO(X)$ be a \lnob  of a subspace
$Y\subseteq X$. Suppose that     ${\rho}(B) \ge \lambda$ for every $B\in {\mc  B}$.
Then $Y$ admits a \nob ${\mc  B}'\subs \RO(X)$.
\end{mlemma}
\begin{proof}
By Lemma~\ref{lem1}(c) and Corollary~\ref{cor}, it is sufficient to show that $\wt Y=Y\sm I(Y)$ has a \nob, where $I(Y)$ denotes the set of isolated points of the subspace $Y$.
We may assume that for every $B\in{\mc  B}$, either $B\cap\wt Y\ne\0$, or $\abs{B\cap Y}=1$ .
Define 
\[
 {\mc  A} =\set{A\in \RO(X): \abs{A\cap Y}\ge\omega},\]
\[ 
\wt{{\mc  B}} = {\mc  B}\cap{\mc  A} = \set{B\in{\mc  B}: B\cap\wt Y\ne\0}.
\]
Then $\wt{\mc  B}$ is an outer base of $\wt Y$.

Now, let  $f:{\mc  B} \to \lambda$ be a function  witnessing that $\mc B$ is \lno.
For each $B\in \mc B$,  that  we can choose a family of regular open subsets
\[
 \{ G^B_\xi:\xi\in\lambda\}\subs\RO(B)
\]
which forms an antichain, i.e. the elements are pairwise incomparable with respect to inclusion.
This is possible because for every $B\in \RO(X)$, we have   $\abs{\RO(B)}=\rho(B)\ge\lambda$, and since $\RO(B)$ is a complete Boolean algebra,  
it follows from   \cite{BaFr82}*{Theorem A} that  $\RO(B)$ contains an antichain of cardinality $|\RO(B)|$.

A  family ${\mc  E} \subset {\mc  P}(X)$ is said to be {\em strongly disjoint} iff $\overline{E_0}\cap\overline{E_1}=\emptyset$ for every pair $\{E_0,E_1\}\in {[{\mc  E}]}^{2}$.

We define 
\[
 {\mc  S} = \set{s\in {}^\omega{\mc  B}:
 \text{$\overline{s(n)}\cap \overline{s(m)}=\empt$  
 for each $n\ne m<{\omega}$ }},\]
 and fix an enumeration
 $
\mc S= \set{s_\alpha: \alpha\in \abs{{\mc  S}}}
$.
Hence, $\ran(s)$  is strongly disjoint for each $s\in \mc S$.

For every $A\in {\mc  A}$, there is an $s\in{\mc  S}$ such, that
\[
 \bigcup\set{\overline{s(n)}:n\in \omega}\subs A.
\]
If $A\in \RO(X)$ and $s\in{\mc  S}$, we say that $\ds s \subs^* A$, if 
\[
 \exists n_0\; \forall n\ge n_0\; s(n) \subs A.
\]
If $s\subseteq^* A$, then let
\[
 n(s,A) = \min\set{n_0: \forall n\ge n_0\; s(n)\subs A}.
\]

\noindent
For every $B\in\wt{\mc  B}$, we define $A(B,0),A(B,1)\in \RO(X)$ as follows: 
\[
 \alpha(B) = \min\set{\alpha: s_\alpha \subs^* B},\qquad
 n(B) = n(s_{\alpha(B)},B).
\]
For $i=0,1$, we denote
\[
 H(B,i) = G^{s_{\alpha(B)}(n(B)+i)}_{f(B)}
\]
that is, we choose the $n(B)$-th and $(n(B)+1)$-th member of $\ds s_{\alpha(B)}$, and from each of them, the $f(B)$-th regular open subset from the corresponding antichain.
Let 
\[
 A(B,i) = B\sm \overline{H(B,i)}.
\]

\begin{Claim}
 For every $B\in\wt{\mc  B}$ and $i\in 2$, 
\[
A(B,i)\in{\mc  A} \mbox{ and } \alpha(A(B,i)) = \alpha(B) \mbox{ and }  n(A(B,i))=n(B)+1+i.
\]
\end{Claim}

\begin{proof}

\noindent
Let $s=s_{\alpha(B)}$ and $m\ge n(B)$.  
Then we have:
\begin{displaymath}\tag{\dag}
\text{$s(m)\subs B$, and $s(m)\subs A(B,i)$ iff $m\ne n(B)+i$ 
}
\end{displaymath}
because $H(B,i)\subs s(n(B)+i)$, so $s(n(B)+i)\not\subs A(B,i)$, while all other $s(m)$ with $m\ne n(B)+i$ remain disjoint from $H(B,i)$, and hence are contained in 
$A(B,i)=B\setm \overline{H(B,i)}$.

This shows that $s\subseteq^* A(B,i)$, and  hence  $|A(B,i)\cap Y|={\omega}$. Moreover,    by definition of ${\alpha}(\cdot)$, we have ${\alpha}(A(B,i))\le {\alpha}(B)$.
On the other hand, since $A(B,i)\subs B$, it follows that $\alpha(B)\le\alpha(A(B,i))$.
Therefore, we conclude 
\begin{displaymath}
{\alpha}(A(B,i))={\alpha}(B).
\end{displaymath}
By definition of $n(\cdot)$,
($\dag$) implies that $n(A(B,i))=n(B)+1+i$.
\end{proof}

Now let
\[
 {\mc  B}_i = \set{A(B,i): B\in\wt{\mc  B}}
\]
for $i\in \{0,1\}$.
Note that  $A(B,0)\cup A(B,1)=B$ since $\cl{H(B,0)}\cap\cl{H(B,1)}=\0$. Hence,  ${\mc  B}_0\cup{\mc  B}_1$ forms  an outer base of $\wt Y$. It remains to prove that each ${\mc  B}_i$ is Noetherian. 
Fix  $i\in \{0,1\}$ and suppose
\[
 \<B_k: k\in\omega\>  \subs\wt{\mc  B} 
\]
is a sequence  such that 
\[
 \<A_k=A(B_k,i): k\in\omega\>\subs \mc B_i
\]
is an increasing sequence.  We will show that the sequence stabilizes; that is, it becomes eventually constant.

We denote
\[
 \alpha_k = \min\set{\alpha: s_\alpha \subs^* A_k},\qquad
n^k = n(s_{\alpha_k},A_k).
\]
These definitions are valid since each  $A_k\in{\mc  A}$.
The sequence $\ds \<\alpha_k: k\in\omega\>$ is a decreasing because the $A_k$'s form an  increasing sequence. Therefore, without loss of generality, we may assume that 
there exists an ordinal  $\alpha$ such that  
\begin{displaymath}
   \text{$\alpha_k = \alpha$ for every  $k\in\omega$.}
\end{displaymath}

Similarly, the sequence $\< n_k: k\in\omega\>$ is non-increasing, 
so it is eventually constant. Thus, 
 we may assume  that
there exists some fixed $n\in {\omega}$ such that 
\begin{displaymath}
\text{$n^k = n$ for all $k\in\omega$.}
\end{displaymath}

By the previously established  Claim, it follows that  
\begin{displaymath}
\text{$\alpha(B_k) = \alpha$ 
and $n(B_k)=n-1-i$
for every $k\in\omega$. 
}
\end{displaymath}

Now suppose, towards a contradiction, that 
$$f(B_k)\ne f(B_{k'})$$ for some $k\ne k'$. Then the corresponding sets   $H(B_k,i)$ and $H(B_{k'},i)$ are incomparable. 
Since $H(B_{k},i)\cup H(B_{k'},i)\subs s_{\alpha}(n-1)$ and 
\begin{align*}
   A_k\cap s_{\alpha}(n-1)=&s_{{\alpha}}(n-1)\setm \overline{H(B_k,i)},\\
   A_{k'}\cap s_{\alpha}(n-1)=&s_{{\alpha}}(n-1)\setm \overline{H(B_{k'},i)},
\end{align*}
it follows that  $A_k$ and $A_{k'}$ are also incomparable, contradicting the 
assumption that $\{A_k:k<{\omega}\}$ is an increasing sequence. Therefore, there is ${\xi}<{\lambda}$ such that
\[
 \forall k\in\omega\; f(B_k)={\xi}.
\]
This implies that the  open set $H=G_{\xi}^{s(n-1)}$ satisfies 
\[
 \forall k\in\omega\;, H(B_k,i) = H.
\]
Since $B_k=A_k\cup \overline{H}$, the sequence $\<B_k: k\in\omega\>$ is  increasing. Since the function $f$ witnessed that $\mc B$ is ${\lambda}$-Noetherian, this sequence must  stabilize.  
Therefore, the sequence $\<A_k:k\in\omega\>$ also stabilizes.
This completes the proof.
\end{proof}

We say that a space $X$  is {\em $\varrho$-homogeneous} iff ${\rho}(B)={\rho}(X)$
for every $\0\ne B\in \RO(X)$.

\begin{theorem}\label{RoHom}
 Let  $X$ be  a $\varrho$-homogeneous $T_3$ space.
 Then $X$ has a \nb\ ${\mc  B^*}\subs \RO(X)$.
\end{theorem}
\begin{proof}
 Let $\lambda = {\rho}(X)$ and $ {\mc  B} = \RO(X)$. 
 Since 
$X$ is 
${\rho}$-homogeneous and 
$\RO(X)$
 is a base for 
$X$, we may apply the Main Lemma for $X$  and $\mc B$  to conclude that 
$X$ has a Noetherian base  $\mc B^*$ contained in 
$\RO(X)$.
\end{proof}

\begin{definition}\label{df:acal}
   Let $X$ be a
   regular topological space $X$.  Let ${\mc  A}(X)$ be a maximal strongly disjoint  family of non-empty, regular-open, ${\rho}$-homogeneous 
 subsets of $X$.  Define the function 
$$
\tr_X:\RO(X)\to \mathcal P({\mc  A}(X))$$
by the formula 
\begin{displaymath}
\tr_{X}(H)=\{A\in {\mc  A(X)}:A\cap H\ne \emptyset\}.
\end{displaymath}
\end{definition}

\begin{lemma}\label{lm:uj}
     Let $X$ be a regular topological space.  
        Then,  
for each $H\in \RO(X)$,
\begin{align}\label{eq:elso}
    \rho(H)=&\prod\{\rho(A): A\in \tr_X(H)\},
    \\\label{eq:masodik}
    \rho(H)\ge& 2^{|\tr_X(H)|},
\end{align} and 
\begin{equation}\label{eq:Bcount}
|\{B\in \RO(X):\tr_X(B)\subs \tr_X(H)\}|\le \rho(H). 
\end{equation}
\end{lemma}

\begin{proof}
 Write $\mc A=\mc A(X)$.  
   For $H_0,H_1\in \RO(X)$ we have 
$H_0=H_1$ iff $H_0\cap \bigcup {\mc  A}=H_1\cap \bigcup {\mc  A}$ because $\bigcup {\mc  A}$ is dense in $X$. 
Thus 
\begin{displaymath}
    \rho(H)=\prod\{\rho(A\cap H): A\in \tr_X(H)\}.
\end{displaymath}
Since every $A\in {\mc  A} $ is ${\rho}$-homogeneous, we have ${\rho}(A)={\rho}(A\cap H)$. 
So we proved  \eqref{eq:elso}.
Since ${\rho}(A)>1$ for $A\in \tr_X(H)$, 
\eqref{eq:elso}
implies  \eqref{eq:masodik}.

To show \eqref{eq:Bcount} let $G=\inte \overline{\bigcup \tr_X(H)}$.
Then ${\rho}(G)={\rho}(H)$ by \eqref{eq:elso}. 

Moreover, if $\tr_X(B)\subs \tr_X(H)$ for some $B\in \RO(X)$, then $B\in \RO(G)$.
Thus \eqref{eq:Bcount} holds. 
\end{proof}

Given a cardinal $\lambda$, define the {\em upper logarithm} of $\lambda$ as follows:
\begin{displaymath}
   \lu(\lambda) = \min\set{\nu : 2^\nu > \lambda}.
\end{displaymath}
Clearly, $2^{<\lu({\lambda})}\le {\lambda}$.

\begin{theorem}\label{acalprops}
 Let $X$ be a $T_3$ space, $\lambda$ a cardinal and define
\[
 X_\lambda = \set{p\in X: \varrho(p,X)=\lambda}.
\]
Then $X_\lambda$ has a \nob ${\mc  B}_{\lambda}$.
\end{theorem}
\begin{proof}
    Define
    \[
     {\mc  B}= \set{B\in \RO(X):{\rho}(B)=\lambda\text{ and } B\cap X_\lambda\ne\0}.
    \]
    Then ${\mc  B}$ is an \ob of $X_\lambda$. 
    To apply the Main Lemma for $X_\lambda$ and ${\mc  B}$,  we must show that 
    ${\mc  B}$ is \lno. 
    
Consider the function $\tr_X$ from  Definition \ref{df:acal}
and let $\tr_{\mc B}=\tr_X\restriction \mc B$.

By the inequality Lemma \ref{lm:uj}.\eqref{eq:masodik},
$2^{\tr_{\mc B}(B)}\le {\lambda}$,
and so 
$|\tr_{\mc B}(B)|<\lu({\lambda})$ for each 
$B\in {\mc  B}$.
Hence $\tr_{\mc B}:\mc B\to {[\mc A]}^{<\lu({\lambda})}$, and 
$\tr_{\mc B}$ is clearly  monotonic.
As we recalled in  Lemma \ref{LuNoether}, the poset  $[{\mc  A}]^{<\lu(\lambda)}$ is $2^{<\lu({\lambda})}$-Noetherian. Since  ${\lambda}\ge 2^{<\lu({\lambda})}$, 
\begin{displaymath}\tag{$*$}
\text{$[{\mc  A}]^{<\lu(\lambda)}$ is ${\lambda}$-Noetherian}.
\end{displaymath}
By the inequality Lemma \ref{lm:uj}\eqref{eq:Bcount}, 
for every $B\in{\mc  B}$ we have $$|\tr_{\mc B}^{-1}\{\tr_{\mc B}(B)\}|\le {\rho}(B)={\lambda},$$ and so  
\begin{displaymath}\tag{$**$}
\text{$\tr_{\mc B}^{-1}\{\tr_{\mc B}(B)\}$
is  \lno for each $B\in \mc B$.} 
\end{displaymath}
By lemma
\ref{MonNoether},
($*$) and  ($**$) together imply that  
  ${\mc  B}$ is \lno, as well.     

Thus,  the Main Lemma applies, and we conclude that  
$X_{\lambda}$ has a Noetherian outer base ${\mc  B}_{\lambda}$.
\end{proof}

\subsection*{Locally small spaces}

Instead Theorem \ref{tm:noether-base} we prove the following stronger result. 

\begin{theorem}\label{tm:noether-stronger}
   Let $\<X,\tau\>$ be a $T_3$ space. Suppose that ${\mu}$ is a cardinal    such that
\[
 \forall p\in X\; \exists G\in\tau\; (p\in G \mbox{ and } \abs G \le \mu).
\]
If there is no strongly inaccessible cardinal less than or equal to  ${\mu}$,
then $X$ has a \nb. 
\end{theorem}
\begin{proof}
    Consider the function $\tr_X$ we defined in Definition \ref{df:acal}.
For every $p\in X$, define 
\[
 \tr_{X}(p) = \min\set{\abs{\tr_X(B)}: p\in B\in \RO(X)},\]
 and let 
 \[\qquad
 \tr = \sup\{\tr_X(p): p\in X\}.
\]
Then clearly $\tr\le \mu$.
Now, for each cardinal  $ \kappa \le \tr$, define
\[
 Z_\kappa = \{{p\in X: \tr_X(p)=\kappa}\}.
\]
\begin{lemma}\label{lm:Z_kappa}
   For every ${\kappa}\le \tr$, the set $Z_\kappa$ has a \nb.
\end{lemma}
\begin{proof}[Proof of the Lemma]
For each ${\lambda}\le {\rho}(X)$, define
\[
 Z_{\kappa,\lambda}=Z_{\kappa}\cap X_{\lambda} = \set{p\in Z_\kappa: \varrho(p,X) = \lambda}.
\]
By Theorem \ref{acalprops}, the subspace $X_{\lambda}$ has a \nobs\ ${\mc  B}_{\lambda}$.
Let 
\begin{displaymath}
{\mc  B}_{{\kappa},{\lambda}}=\{B\in {\mc  B}_{\lambda}: 
\abs{\tr_X(B)}=\kappa \mbox{ and } {\rho}(B) = \lambda\}.
\end{displaymath}
Then ${\mc  B}_{{\kappa}, {\lambda}}$ is a \nob of $Z_{{\kappa},{\lambda}}$.

Let
\[
 {\mc  B} = \cup\set{{\mc  B}_{{\kappa},\lambda}: Z_{\kappa,\lambda} \ne \0}.
\]
Then ${\mc  B}$ is a \ob for $Z_\kappa$. 

\begin{Claim}
   ${\mc  B}$ is $2^\kappa$-Noetherian.
\end{Claim}

\begin{proof}[Proof of the Claim]

Take an arbitrary $B\in{\mc  B}$. Then there exists  a unique $\lambda$ such that $B\in{\mc  B}_\lambda$. 
If $\tr_X(B')=\tr_X(B)$ for some $B'\in {\mc  B} $, then ${\rho}(B)={\rho}(B')$ by Lemma \ref{lm:uj}(1), so 
$B'\in {\mc  B}_\lambda$. Thus,  for every $B\in{\mc  B}$, we have  
\[
 \tr_X^{-1}(\tr_X(B))\subs {\mc  B}_\lambda,
\]
which  implies that $\tr_X^{-1}(\tr_X(B))$ is  Noetherian.

The function $\tr_X:{\mc  B}\to [{\mc  A}(X)]^\kappa$ is monotonic, and the poset  
$[{\mc  A}(X)]^\kappa$
is $2^\kappa$-Noetherian by Lemma \ref{LuNoether}. 
 Thus, by Lemma \ref{MonNoether}, ${\mc  B}$ is a 
 $2^\kappa$-Noetherian, as well. So we proved the Claim.
\end{proof}
Furthermore, for every $B\in{\mc  B}$
\begin{displaymath}
   \rho(B) \ge 2^{|\tr_X(B)|}=2^\kappa
\end{displaymath}
by Lemma \ref{lm:uj}\eqref{eq:masodik}. Since $\mc B$ is $2^{\kappa}$-Noetherian by the Claim, 
we can apply the Main Lemma to $Z_\kappa$ and $\mc B$, and conclude that $Z_\kappa$ has a \nob. This completes the proof  of the Lemma. 
\end{proof}

\noindent
For a cardinal $\kappa$, let us denote by $\ds Z_{< \kappa} = \{p\in X: \tr_X(p)<\kappa\}$. Then clearly
 $X=Z_{<r}\cup Z_r$. By the Lemma \ref{lm:Z_kappa} above, $Z_r$ has a \nob. Therefore, 
 by Lemma \ref{lem1}(c), it suffices to prove that for every $\kappa\le \tr$ 
 \begin{displaymath}\tag{$\star_{\kappa}$}
 \text{the subspace $Z_{<\kappa}$ also has a \nob.}
 \end{displaymath}
 We prove this by 
transfinite induction of $\kappa$. 

\noindent
{\bf Case 1.} $\kappa=\omega$.

If $\tr_X(p)<{\omega}$, then $\tr_X(p)=1$ because $\mc A(X)$ is strongly disjoint. 
Moreover, if $B\in RO(X)$  and  $\tr_X(B)=\{A\}$ for some $A\in \mc A_X$,
then $B\setm \overline{A}=\empt$, and so $B\subs \inte \overline{A}=A$.
Thus $Z_{<{\omega}}=Z_1\subs \bigcup \mc A$. 

On the other hand, $\tr_X(p)=1$ for each $p\in A\in \mc A_X$.
Hence $Z_{<\omega}=Z_1=\cup\mathcal A$, 
 that is $Z_{<\omega}$ is a disjoint union of $\varrho$-homogeneous open sets. 
 Using Theorem \ref{RoHom} and Lemma \ref{NoUnion},
we conclude that  $Z_{<{\omega}}$ has a \nob.

\noindent
{\bf Case 2.} $\kappa=\mu^+$ for some $\mu\ge \omega$ cardinal.

In this case,  $Z_{<\kappa} = Z_{<\mu}\cup Z_\mu$.
By the induction hypothesis ($\star_{\mu}$)  and by Lemma \ref{lm:Z_kappa}, both $Z_{<\mu}$ and $Z_{\mu}$ have \nob{s}, so  $Z_{<{\kappa}}$ does as well, by Lemma \ref{lem1}(c).

\noindent
{\bf Case 3.} $\omega<\kappa\le \tr$ and $\kappa$ is a limit cardinal.

The cardinal   $\kappa$ is not strongly inaccessible  because ${\kappa}\le {\mu}$. Hence,  
there exists  a cardinal $\nu < \kappa$ such that $\cf \kappa \le 2^\nu$.
Let
\[
 \<\kappa_\xi: \xi<\cf\kappa\>
\]
be a closed cofinal  sequences of cardinals in $\kappa$, with $\kappa_0 = \nu$.
Then 
\[
 Z_{<\kappa} = Z_{<\kappa_0}\cup \bigcup\set{Y_\xi: \xi < \cf\kappa},
\]
where
\[
 Y_\xi = \bigcup\set{Z_\sigma: \kappa_\xi\le\sigma<\kappa_{\xi+1}}\subs Z_{<\kappa_{\xi+1}}.
\]
Let $$Y=\bigcup\set{Y_\xi: \xi < \cf\kappa}.$$ By the induction hypothesis, $Z_{<\kappa_0}$ has a \nob, so it suffices  to show, that $Y$ has one too.
For every $\xi<\cf \kappa\;$,  the subspace $Z_{<\kappa_{\xi+1}}$ has a \nob by the induction hypothesis, so its   subset $Y_\xi$ does as well. 

Hence, the union 
\begin{displaymath}
\mc B=\bigcup_{{\xi}<\cf({\kappa})}\mc B^{\xi},
\end{displaymath}
where each 
$\mc B^{\xi}$  is a Noetherian outer base for 
$Y_{\xi}$, is a $\cf({\kappa})$-Noetherian outer base for $Y$. 

Since $\cf ({\kappa})\le 2^{{\kappa}_0}$, this outer base is also $2^{{\kappa}_0}$-Noetherian.

Moreover, for every 
$B\in{\mc  B}$, by the definitions of $Z_\sigma$-s,
\[
 \rho(B) \ge 2^{|\tr_X(B)|}\ge 2^{\kappa_0}.
\]
by Lemma \ref{lm:uj}\eqref{eq:masodik}.
Therefore,  we can apply the Main Lemma 
for $Y$ and $\mc B$
to conclude that  $Y$ has a \nob. This completes the induction and hence  the proof of the theorem.
\end{proof}


\section{Noetherian type}\label{sc:type}
The $G_{\delta}$-modifications of  topological spaces  have been 
extensively studied in the literature. 
A natural problem in this area is 
to establish an upper  bound for a given cardinal invariant on 
$X_{\delta}$ in terms of its value on $X$ (see e.g. 
\cites{BeSp19,Ar06,Us19,AuBe14,DoJuSoSzWe19,BeSp21}). 
Milovich and Spadaro (\cite{KoMiSp14})  considered  this problem for the Noetherian type.
\begin{theorem}[Spadaro, see \cite{KoMiSp14}]\label{spadaro} If GCH holds and  $X$ is  a compact space such that $\Nt(X)$ has
   uncountable cofinality, then   $\Nt(\deltop X) \le 2^{N t(X)}$.  
   \end{theorem}
   Since Milovich, \cite{KoMiSp14}, proved that if  $X$  is  compact 
   dyadic homogeneous space then $\Nt(X)=\omega$,
    Theorem \ref{spadaro} does not  apply  to the Cantor cubes.
   However, Milovich  proved that $\Nt ((D(2)^{\kappa}{})_{\delta})=\omega_1$ for $\kappa<\alo$
   under GCH.
   So the  simplest unsettled case remained $\kappa=\alo$. 
This problem was also addressed in \cite{KoMiSp14}.   

\begin{definition}\label{df:noetherian-type-nt}
   Write {\bf (Nt)} iff $\Nt((D(2)^{\alo})_{\delta})= {\omega}_1$.
   \end{definition}

\begin{theorem}[\cite{KoMiSp14}*{Corollary 3.24}]\label{tm:Mi}
 GCH +
$\square_{\aleph_{\omega}}$ implies {\bf Nt}. 
\end{theorem}
On the other hand, in \cite{soukup2010notenoetheriantypespaces} the first author of the present paper proved: 
\begin{theorem}[\cite{soukup2010notenoetheriantypespaces}]\label{tm:SL-1}
If GCH +  $(\aleph_{{\omega}+1},\aleph_{{\omega}})\doublearrow  (\aleph_1,\aleph_0)$
hold, then 
{\bf Nt} fails.
\end{theorem}

This phenomenon is not exceptional. 
In recent decades, several topological and set-theoretical statements have been shown to be 
independent of ZFC by demonstrating that they hold under GCH + $\square_{\aleph_{\omega}}$, 
while their negations follow from 
GCH and a suitable version of the Chang Conjecture, --- specifically, the principle $(\aleph_{{\omega}+1},\aleph_{{\omega}})\doublearrow  (\aleph_1,\aleph_0)$.
In what follows,  we will formulate several  such statements and  investigate the possible  
implications 
between these statements. 
To do so,  we begin by  recalling some relevant notions and notations.

\begin{definition}[\cite{JuShSo88}]\label{df:splendid}
    A regular topological 
    space $X$ is called {\em splendid} if it is countably compact, locally
    countable, and  
    {\em $\omega$-fair}, i.e. $|\overline{Y}|={\omega}$  for each 
    $Y\in {[X]}^{{\omega}}$. 

    Write {\bf (SPL)} if 
there is a splendid space of size $\aleph_{{\omega}+1}$.
\end{definition}

\begin{definition}[\cite{FuKoSh96}] \label{df:wfn}
   A poset $\<P,\le\>$    
    has the {\em weak Freese-Nation} property iff there is a function $F:\mc P\to [{P}]^{{\omega}}$ such that for each 
   $p,q\in P$ with $p\le q$ there is $r\in F(p)\cap F(q)$
   with $p\le r\le q$.
   
   Write {\bf (wFN)} iff the poset 
   $\<{[\aleph_{\omega}]}^{{\omega}},\subseteq\>$
   has the weak Freese-Nation property. 
   \end{definition}

\begin{definition}[\cite{GoJuSh91}]\label{df:saturated}
An almost disjoint family  $\mc A\subs {[{\kappa}]}^{{\omega}}$ is called {\em saturated} iff for each $X\in {[{\kappa}]}^{{\omega}}$
either there is $A\in \mc A$ with $A\subs X$, or 
$X\subset^* \bigcup \mc A'$ for some finite $\mc A'\subs \mc A$.

Write {\bf (SAT)} iff there is a saturated family 
$\mc A\subs {[\aleph_{\omega}]}^{{\omega}}$.
\end{definition}

\begin{definition}[\cite{Ne88}]\label{df:perm_group}

   A permutation group  $G$ on an uncountable  set $A$ is  {\em ${\omega}$-homogeneous}
   iff for all  $X,Y\in\br A;{\omega};$ there is
   a $g\in G$ with $g''X=Y$.
   $G$ is 
   {\em ${\omega}$-transitive}
   iff for any countable injective function  $f$  with $\dom(f)\cup\ran(f)\in 
   \br A; {\omega};$  
    there is
   a $g\in G$ with $f\subs g$.
   
Write {\bf (HnT)} if there is an ${\omega}$-homogeneous, but not 
${\omega}$-transitive permutation group on $\aleph_{\omega}$.  
\end{definition}

      \begin{definition}[\cite{KoMiSp14}]\label{df:sparse}
        A family $\mc A\subs {[{\kappa}]}^{{\omega}}$ is 
        {\em $({\omega}_1,{\omega}_1)$-sparse} if
        $|\bigcup \mc H|\ge {\omega}_1 $ for each $\mc H\in {[\mc A]}^{{\omega}_1}$.

        Write {\bf (SPA)} if there is a cofinal, $({\omega}_1,{\omega}_1)$-sparse
        family $\mc A\subs {[\aleph_{\omega}]}^{{\omega}}$.
    \end{definition}

    \begin{definition}\label{df:;small}
      A family $\mc A\subs {[{\kappa}]}^{{\omega}}$ is 
      {\em locally small} if
      $
      \{X\cap A:A\in \mc A\}$ is countable  for each $X\in {[{\kappa}]}^{{\omega}}$.

      Write {\bf (CLS)} if there is a cofinal, locally small   family $\mc A\subs {[\aleph_{\omega}]}^{{\omega}}$.
  \end{definition}
Clearly, a locally small family is $({\omega}_1,{\omega}_1)$-sparse. So CLS implies SPA.

      The following list summarizes our previous knowledge: 
\begin{theorem}[GCH]\label{tm:summ}
   \makebox[1pt]{}
\begin{enumerate}[(A)]
\item {\em If  $\square_{\aleph_{\omega}}$ holds,
then \smallskip 
\begin{enumerate}[(a)]
\item SPL holds by  \cite{JuNaWe79}*{Theorem 11}, \smallskip 
\item Nt holds   by \cite{KoMiSp14}*{Corollary 3.24}, \smallskip
\item wFN holds by \cite{FuSo97}*{Corollary 11},  \smallskip
\item HnT holds by \cite{ShSo23}*{Theorem 2.5}, \smallskip
\item SAT holds by \cite{GoJuSh91}*{Conclusion}.
\end{enumerate}} \smallskip
\item {\em SPA  iff  Nt } by
\cite{soukup2010notenoetheriantypespaces} and \cite{KoMiSp14}. \smallskip
\item {\em If  $(\aleph_{{\omega}+1},\aleph_{{\omega}})\doublearrow  (\aleph_1,\aleph_0)$ holds, then \smallskip
\begin{enumerate}[(a)]
\item SPL fails by \cite{JuShSo88}*{Corollary 2.2}, \smallskip
\item  SPA, and  hence Nt as well, fails by   
\cite{soukup2010notenoetheriantypespaces},  \smallskip
\item  wFn fails by 
\cite{FuSo97}*{Theorem 12}, and \smallskip
\item SAT may hold.
\end{enumerate} } \smallskip 
\end{enumerate}
\end{theorem} 

To verify  (C)(d), it is enough to observe that, by \cite{HaJuSo87}, 
if one iteratively adds ${\omega}_1$ many dominating reals to a ground model, then in the resulting generic extension,
 (SAT) holds. Moreover, the partition relation $(\aleph_{\omega+1}, \aleph_{\omega}) \doublearrow (\aleph_1, \aleph_0)$ is preserved by any c.c.c.\ forcing.

As the main  result of this section, 
 we establish the following implications among these principles,  assuming GCH.

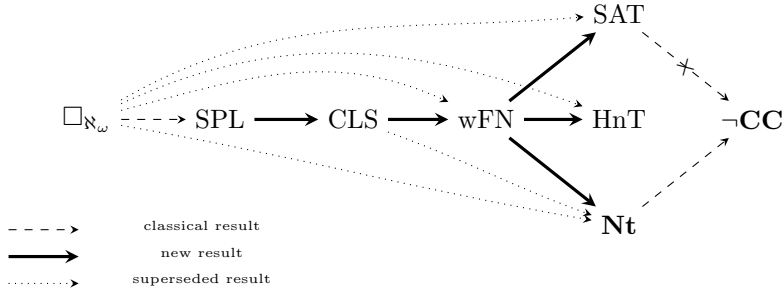
\begin{figure}[ht]

\begin{tikzpicture}[scale=0.7,
   old/.style={-stealth,dashed},
   oldeq/.style={stealth-stealth,very thin},
   new/.style={-stealth,very thick},
   avitt/.style={-stealth,dotted},
   ]
    \node (box)  at (0,0)  {  $\square_{\aleph_{\omega}}$}; 
    \node (spl)  at (2.5,0) {SPL}; 
    \node (cls) at (5,0)  {CLS} ; 
    \node (wfn)  at (7.5,0)   {wFN}; 
    \node (hnt)  at (10,0)  {HnT}; 
    \node (ncc) at (12.5,0) {$\neg$\bf CC};
    \node (sat) at (10,2) {SAT};
    \node (nt) at (10,-2) {\bf Nt};

    \path[new] (spl) edge (cls) ; 
    \path[new] (cls) edge (wfn) ; 
   \path[new] (wfn) edge (nt) ; 
    \path[new] (wfn) edge (hnt); 
   \path[new] (wfn) edge (sat) ; 

   \path[old] (box) edge (spl) ; 
   \path[old] (nt) edge (ncc) ; 
    \path[old] (sat) edge  node[rotate=30,dotted] {$\times$}     (ncc) ; 
    \draw[avitt]  (box)  edge[out=15, in=155]  (wfn);
    \draw[avitt]  (box) edge[out=30, in=185] (sat); 
  \draw[avitt]  (box) edge[out=350, in=170] (nt); 
  \path[avitt] (cls) edge[out=340, in=160] (nt) ; 
  \path[avitt] (box) edge[out=25, in=160] (hnt) ; 
  \matrix[above right,matrix of  nodes,nodes in empty cells,row sep
   = 0em,
   column sep = 1em] (m) at (-2,-3.5) {
   &&&\tiny classical result\\
   &&&\tiny new result\\
   &&&\tiny superseded result \\  };
   \draw[old]  (m-1-1) -- (m-1-3);
   \draw[new]  (m-2-1) -- (m-2-3);
   \draw[dotted,-stealth]  (m-3-1) -- (m-3-3);
\end{tikzpicture} 
\caption{Implication diagram of combinatorial, topological and algebraic statements  at $\aleph_{\omega}$ and $\aleph_{{\omega}+1}$. }\label{fig:1}
\end{figure}

\begin{theorem}\label{tm:main-implication}
Assume GCH. Then
\begin{enumerate}[(1)]
\item \label{cs:spl2cls} SPL implies CLS.
\item \label{cs:cls2wfn}  CLS implies wFn.
\item\label{cs:wfn2sat} wFN implies SAT.
\item \label{cs:wfn2nt}wFN  implies Nt. 
\item \label{cs:wfn2hnt} wFN  implies HnT. 
\end{enumerate}
\end{theorem}
Figure \ref{fig:1} summarizes both the previously known implications and the new results.  
We use     {\bf $\bm\neg$CC} to denote the failure of   
    the Chang conjecture    
    $(\aleph_{{\omega}+1},\aleph_{{\omega}})\doublearrow  (\aleph_1,\aleph_0)$.

\begin{proof}
(\ref{cs:spl2cls}). 
We need the following statement.       
\begin{lemma}   \label{lm:splendid2nice}
   If $X$ is a splendid space, then the family $\mc U$   
   of compact open 
   subsets of $X$ is cofinal in ${[X]}^{{\omega}}$ and locally small.  		
   \end{lemma}

This  statement follows from the proof of \cite{ShSo23}*{Lemma 2.4}, 
but  we include it here  for completeness. 

   \begin{proof} Since $X$ is locally countable, $\mc U\subs {[X]}^{{\omega}}$.
      
       Let $A\in [X]^{\omega}$. Then $\overline A$ is countable, and hence compact.
        Since a splendid space is zero-dimensional, $A$ can be covered by 
       finitely many compact open set, and so $A$ can be covered by 
       an element of $\mc U$. Thus ${\mc U}$ is  cofinal in 
       $\<[X]^{\omega},\subseteq\>$.  
       
       To verify that $\mc U$ is 
        locally small, observe that every   $U\in \mc U$ is a compact countable space, 
       hence homeomorphic to a countable successor ordinal.  Therefore,  $U$
       has only countably many compact open subsets, so
        $\mc U\cap \mc P(U)$ is countable. Hence $\mc U$ is 
        $({\omega}_1,{\omega}_1)$-sparse.
        Since $\mc U$ is closed under intersection, it follows that $\mc U$ is locally small as well. 
       \end{proof}

(SPL) and Lemma \ref{lm:splendid2nice} together yield a  family $\mc B$ 
which is locally small and cofinal in $[\aleph_{{\omega}+1}]^{\omega}.$
Then the family $\mc A=\{B\cap \aleph_{\omega}:B\in \mc B\}$ is also locally small, and clearly cofinal in 
${[\aleph_{\omega}]}^{{\omega}}$.

\noindent (\ref{cs:cls2wfn}). 
Let $\mc A=\{A_{\alpha}:{\alpha}<\aleph_{{\omega}+1}\}
\subs {[\aleph_{{\omega}}]}^{{\omega}}$ be a cofinal, locally small family.

Since
$\mc A$ is locally small, 
for each ${\alpha}$ we can pick $I_{\alpha}\in {[{\alpha}]}^{{\omega}}$
such that 
\begin{displaymath}
\{A_{\alpha}\cap A_{\eta}:{\eta}<{\alpha}\}=
\{A_{\alpha}\cap A_{{\zeta}}:{\zeta}\in I_{\alpha}\}. 
\end{displaymath}

For each ${\alpha}<\aleph_{{\omega}+1}$
fix an enumeration $\{a_{{\alpha},i}:i<{\omega}_1\}={[A_{\alpha}]}^{{\omega}}$,
let 
$\mc A_{\alpha}={[A_{\alpha}]}^{{\omega}}\setm \bigcup_{{\xi}<{\alpha}}
{[A_{\xi}]}^{{\omega}}$
and 
$\mc A_{\le{\alpha}}=\bigcup_{{\zeta}\le{{\alpha}}}\mc A_{\zeta}=
\bigcup_{{\zeta}\le{{\alpha}}}[A_{\zeta}]^{\omega}$.

By transfinite induction on ${\alpha}$,
define functions $F_{\alpha}: \mc A_{\alpha}\to [\mc A_{\le{\alpha}}]^{\omega} $ for ${\alpha}<\aleph_{{\omega}+1}$
as follows.  

To simplify the construction 
declare that $F_{{\alpha}}(x)=\empt$ for each $x\in {[A_{\alpha}]}^{<{\omega}}$.  
Write  $F_{<{\alpha}}=\bigcup_{{\xi}<{\alpha}} F_{\xi}$,
and for $a_{{\alpha},j}\in \mc A_{\alpha}$
let 
\begin{displaymath}
F_{\alpha}(a_{{\alpha},j})=
\{a_{{\alpha},i}:i\le j\}\cup\bigcup_{{\zeta}\in I_{\alpha}}
F_{<{\alpha}}(a_{{\alpha},j}\cap A_{\zeta}).
\end{displaymath}
Since $a_{{\alpha},j}\cap A_{\eta}\in  {[A_{\eta}]}^{\le{\omega}}\subs \dom(F_{<{\alpha}})$, the recursive definition above is meaningful.

To verify that the function $F=\bigcup\{F_{\alpha}:{\alpha}<\aleph_{{\omega}+1}\}$  witnesses (wFN), we show,
by transfinite induction on ${\alpha}$,
 that
\begin{enumerate}[($\star_{\alpha}$)]
\item if $a_{{\alpha},j}\in \mc A_{\alpha}$, $a_{{\eta},i}\in \mc A_{\eta}$,  
and $a_{{\eta},i}\subs a_{{\alpha},j}$ then\\ \makebox[1cm]{}
there is $z\in F(a_{{\alpha},j})\cap F(a_{{\eta},i}) $ with 
$a_{{\eta},i}\subseteq z\subseteq a_{{\alpha},j}$.
\end{enumerate}
Since $\mc A_{\le {\alpha}}$
is closed under subsets, we have ${\eta}\le {\alpha}$.   

If ${\eta}={\alpha}$, then  
$a_{{\alpha},\min(i,j)}\in
F_{\alpha}(a_{{\alpha},i})\cap F_{\alpha}({a_{{\alpha},j}})=
F(a_{{\alpha},i})\cap F({a_{{\alpha},j}})
$ by the definition of $F_{\alpha}$.

So we can assume that ${\eta}<{\alpha}$. 
Pick ${\zeta}\in I_{\alpha}$ such that 
$A_{\alpha}\cap A_{\eta}=A_{\alpha}\cap A_{\zeta}$.

Then $a_{{\eta},i}\subsq a_{{\alpha},j}\cap A_{\eta}=
a_{{\alpha},j}\cap A_{\zeta}
\in \mc A_{\le {\zeta}}$, and so applying the inductive
assumption for $a_{{\eta},i}$ and $a_{{\alpha},j}\cap A_{{\zeta}}$ 
there is $z\in F(a_{{\eta},i})\cap 
F(a_{{\alpha},j}\cap A_{\zeta})$
with $a_{{\eta},i}\subsq z\subsq a_{{\alpha},j}\cap A_{\zeta}$.  
By the definition of $F_{{\alpha}}(a_{{\alpha},j})$, 
we have $z\in F_{<\alpha}(a_{{\alpha},j}\cap A_{\zeta})\subsq 
F_{\alpha}(a_{{\alpha},j}).$
Thus $z\in F(a_{{\eta},i})\cap F(a_{{\alpha},j})$ and 
$a_{{\eta},i}\subsq z\subsq a_{{\alpha},j}$.

Hence, we carried out the inductive step, so $(\star_{\alpha})$ holds for all 
${\alpha}<\aleph_{{\omega}+1}$. Thus $F$ is a wFN-function, so 
(2) is proved. 

\medskip
To prove  (\ref{cs:wfn2sat}), (\ref{cs:wfn2nt}) and (\ref{cs:wfn2hnt})
we need the following technical lemma: 

\begin{lemma}[GCH]\label{lm:wFN2enu}
Given any function   $F: {[\aleph_{{\omega}}]}^{{\omega}}\to [{[{\aleph_{\omega}}]}^{{\omega}}]^{\omega}$, 
there is an enumeration 
$\{a_{{\alpha},j}:{\alpha}<\aleph_{{\omega}+1},j<{\omega}_1\}$ 
of ${[\aleph_{\omega}]}^{{\omega}}$ without repetition such that 
\begin{enumerate}[(\dag)]
\item if $x \in [a_{{\alpha},j}]^{\omega}$
then $F(x)\subs \{a_{{\eta},i}:{\eta}\le {\alpha}, i<{\omega}_1\}$
for each ${\alpha}<\aleph_{{\omega}+1}$ and $j<{\omega}$.
\end{enumerate}
\end{lemma}

\begin{proof}
   We say that $A\subs \aleph_{{\omega}+1}$ is {\em $F$-closed} 
   iff $F(x)\subs {[A]}^{{\omega}}$ for each $x\in {[A]}^{{\omega}}$.
   Since GCH holds, 
   \begin{displaymath}
      \mc A= \{A\in {[\aleph_{\omega}]}^{{\omega}_1}:
      \text{$A$ is $F$-closed}\}
      \end{displaymath}
   is cofinal in $[\aleph_{\omega}]^{{\omega}_1}$,
   and 
   we can choose  
     $\{A_{\alpha}:{\alpha}<\aleph_{{\omega}+1}\}$
   as  a cofinal subfamily  of $\mc A$
   such that 
   for each ${\alpha}$ the family 
   \begin{displaymath}
   \mc A_{\alpha}={[A_{\alpha}]}^{{\omega}}\setm \bigcup_{{\gamma}<{\alpha}}
   {[A_{\gamma}]}^{{\omega}}
   \end{displaymath}
   is uncountable. 
   Let  $\{a_{{\alpha},j}:j<{\omega}_1\}$ be an enumeration 
 of  $\mc A_{\alpha}$  without repetition for 
   ${\alpha}<\aleph_{{\omega}+1}$.
Then the  enumeration  $\{a_{{\alpha},j}:{\alpha}<\aleph_{{\omega}+1},j<{\omega}_1\}$ 
satisfies the requirements. 
\end{proof}

Using  Lemma~\ref{lm:wFN2enu} we first prove (\ref{cs:wfn2sat}), then  (\ref{cs:wfn2nt}) and (\ref{cs:wfn2hnt}).

\medskip
\noindent  (\ref{cs:wfn2sat}).
For families  $\mc H,\mc A\subs {[\aleph_{{\omega}}]}^{{\omega}}$
we will say that $\mc H$ is {\em $\mc A$-saturated}
iff for each $a\in \mc A$ either there is $h\in \mc H$
with $h\subs a$, or $a\subs^* \bigcup \mc H'$
for some finite $\mc H'\subs \mc H$.

Fix a wFN function $F: {[{\kappa}]}^{{\omega}}\to [{[{\kappa}]}^{{\omega}}]^{\omega}$, and 
consider the enumeration $\{a_{{\alpha},j}:{\alpha}<\aleph_{{\omega}+1},j<{\omega}_1\}$ 
of ${[\aleph_{\omega}]}^{{\omega}}$ from Lemma \ref{lm:wFN2enu}.

Consider the lexicographical  ordering $\triangleleft$ of 
$\aleph_{{\omega}+1}\times {\omega}_1$. 
By transfinite recursion on $\triangleleft$,
 we define the elements of the set $D\subs \aleph_{{\omega}+1}\times {\omega}_1$ 
 and $h_{{\alpha},j}\in {[a_{{\alpha},j}]}^{{\omega}}$
 for $\<{\alpha},j\>\in D$ such that writing 
$$\mc H_{{\alpha},j}=\{h_{{\eta},i}:\<{\eta},i\>\in D,\<{\eta},i\>\triangleleft \<{\alpha},j\> \}$$
and  $$\mc A_{{\alpha},j}=\{A_{{\eta},i}:\<{\eta},i\>\triangleleft \<{\alpha},j\> \}$$ 
we have that 
\begin{displaymath}\tag{$\star{}_{{\alpha},j}$}
\text{$\mc H_{{\alpha},j}$ is almost disjoint and $\mc A_{{\alpha},j}$-saturated. }
\end{displaymath}

Assume that for each 
$\<{\eta},i\>\triangleleft \<{\alpha},j\>$ we have decided 
whether $\<{\eta},i\>\in D$ and constructed 
$h_{{\eta},i}$ whenever  $\<{\eta},i\>\in D$.

Consider the family  
\begin{displaymath}
\mc Z=\{a_{{\eta},i}\in F(a_{{\alpha},j}): 
\<{\eta},i\>\triangleleft \<{\alpha},j\>, 
a_{{\eta},i}\subs a_{{\alpha},j}\}
\end{displaymath}
We should distinguish cases.

\medskip

\noindent {\bf Case 1.}
{\em There are $z\in \mc Z$ and $h\in \mc H_{{\alpha},j}$  such that $h\subs z$.}

In this case, let 
$\<{\alpha},j\>\notin D$.

\medskip

\noindent {\bf Case 2.}
{\em There are no $z\in \mc Z$ and 
$h\in \mc H_{{\alpha},j}$ 
such that $h\subs z$. 
}

Since $\mc H_{{\alpha},j}$ is 
$\mc A_{{\alpha},j}$-saturated and $\mc Z\subs \mc A_{{\alpha},j}$, it follows that  
for each  $z\in \mc Z$ there exists  a finite subset 
$\mc K_z\subs \mc H_{{\alpha},j}$
such that $z\subs^* K_z$, where  
 $K_z=\bigcup \mc K_z$.
\medskip

\noindent {\bf Case 2.1.}
{\em A finite subset of the countable family 
$$\mc K=\{K_z:z\in \mc Z\}\cup \{h_{{\alpha},i}:
i<j, \<{\alpha},i\>\in D\}$$
covers $a_{{\alpha},j}$ mod finite.}

In this case, let $\<{\alpha},j\>\notin D$.

\medskip

\noindent {\bf Case 2.2.}
{\em No finite subset of the countable family 
$$\mc K=\{K_z:z\in \mc Z\}\cup \{h_{{\alpha},i}:
i<j, \<{\alpha},i\>\in D\}$$
covers $a_{{\alpha},j}$ mod finite.}

In this case,   put $\<{\alpha},j\>\in D$ and choose 
$h_{{\alpha},j}\in {[a_{{\alpha},j}]}^{{\omega}}$ such that 
\begin{displaymath}\tag{$\circ$}
\text{$h_{{\alpha},j}\cap K$ is finite for each 
$K\in \mc K$.}
\end{displaymath}

We should verify that ($\star_{{\alpha},j+1})$ holds.

\begin{Claim}
   $\mc H_{{\alpha},j+1}$ is almost disjoint. 
\end{Claim}

\begin{proof}[Proof of the Claim]

We may assume that $\<{\alpha},j\>\in D$, i.e. 
we are in Case 2.2.

Assume, towards a contradiction,  that 
there exists $h_{{\eta},\ell}\in \mc H_{{\alpha},j}$ such that 
\begin{displaymath}\tag{$\bullet$}
\text{$h_{{\eta},\ell}\cap h_{{\alpha},j}$ is infinite.}
\end{displaymath}
Since $\{h_{{\alpha},i}:i<j, \<{\alpha},i\>\in D\}\subs \mc K$, we have 
 ${\eta}<{\alpha}$ by $(\circ)$. 

By \ref{lm:wFN2enu}(\dag), we have 
${[a_{{\eta},\ell}]}^{{\omega}}
\subs \{a_{{\xi},k}: {\xi}\le {\eta}, k<{\omega}_1\}$,
so $h_{{\eta},\ell}\cap h_{{\alpha},j}=a_{{\xi},k}$ for some ${\xi}\le {\eta}$ and $k<{\omega}_1$.

Now choose $z\in F(a_{{\xi},k})\cap F(a_{{\alpha},j})$
such that  $$a_{{\xi},k}\subsq z\subsq a_{{\alpha},j}.$$

By \ref{lm:wFN2enu}(\dag), we have 
$z\in \{a_{{\zeta},m}: {\zeta}\le {\xi}, m<{\omega}_1\}$, so $z\in \mc Z$.
Since we are in Case 2, $z\subsq K_z$, where $K_z=\bigcup \mc K_z$ for some finite subset 
$\mc K_z$ of $\mc H_{{{\alpha},j}}$.

Since  $h_{{\eta},\ell}\cap h_{{\alpha},j}=a_{{\xi},k}\subs z\in \mc Z$, 
it follows that 
$h_{{\eta},\ell}\cap h_{{\alpha},j}\subs K_z\in \mc K$, and so 
$h_{{\alpha},j}\cap (h_{{\eta},\ell}\cap h_{{\alpha},j})\subs 
h_{{\alpha},j}\cap K_z$ is finite by $(\circ)$, which contradicts ($\bullet$).
\end{proof}

Next we show that $\mc H_{{\alpha},j+1}$ 
is  $\mc A_{{\alpha},j+1}$-saturated. 
By the inductive hypothesis, we should consider only the set $a_{{\alpha},j}$. 
In Case 1, there is $h\in \mc H_{{\alpha},j}$
with $h\subs a_{{\alpha},j}$.

In Case 2.1 
$a_{{\alpha},j}$ is  covered, mod finite,  by finitely many elements of
$\mc K$. Since every element of $\mc K$
is covered by finitely many elements of $\mc H_{{\alpha},j+1}$, we are done.

Finally, 
in Case 2.2 we have $h_{{\alpha},j}\subs a_{{\alpha},j}$.

So we verified that  ($\star_{{\alpha},j+1})$ holds, that is, the inductive 
construction can be carried out, and so 
the family 
\begin{displaymath}
\mc H=\{h_{{\alpha},j}:\<{\alpha},j\>\in D\}
\end{displaymath}
is almost disjoint and ${[\aleph_{\omega}]}^{{\omega}}$-saturated, i.t. it is a saturated subset of  ${[\aleph_{\omega}]}^{{\omega}}$.

\medskip
\noindent (\ref{cs:wfn2nt}) and (\ref{cs:wfn2hnt}).

Fix a wFN function $F: {[{\kappa}]}^{{\omega}}\to [{[{\kappa}]}^{{\omega}}]^{\omega}$, 
and consider the corresponding enumeration $\{a_{{\alpha},j}:{\alpha}<\aleph_{{\omega}+1},j<{\omega}_1\}$ 
of ${[\aleph_{\omega}]}^{{\omega}}$ from Lemma \ref{lm:wFN2enu}.

Let $\triangleleft$ denote the lexicographical 
order of $\aleph_{{\omega}+1}\times {\omega}_1$, and 
define the family $\mc H\subs {[\aleph_{\omega}]}^{{\omega}}$
as follows:
\begin{displaymath}
a_{{\alpha},j}\in \mc H\text{ iff }
\forall \<{\gamma},k\>\triangleleft \<{\alpha},j\> \
(a_{{\alpha},j}\not\subset a_{{\gamma},k}).
\end{displaymath}

\smallskip
\noindent{\bf Claim 1.}{ \em $\mc H$ is cofinal
in ${[\aleph_{\omega}]}^{{\omega}}.$}

Indeed, given any $x\in {[\aleph_{\omega}]}^{{\omega}}$
let 
\begin{displaymath}
   {\gamma}^*=\min\{{\gamma}:x\subs A_{{\gamma}}\}\text{ and }
k^*=\min\{k<{\omega}_1: x\subset a_{{\gamma},k}\}.
\end{displaymath}
Then $x\subs a_{{\gamma}^*,k^*}\in \mc H$.

\smallskip
\noindent{\bf Claim 2.}{ \em $\mc H$ is $({\omega}_1,{\omega}_1)$-sparse.}

By transfinite induction on the well-order $\triangleleft$
we prove that for each $\<{\alpha},j\>\in \aleph_{{\omega}+1}\times {\omega}_1$ the family
\begin{displaymath}
\text{ $\mc H_{{\alpha},j}=\{a_{{\eta},k}\in \mc H: a_{{\eta},k}\subs a_{{\alpha},j}\}$}
\end{displaymath}
is countable. 

Fix $\<{\alpha},j\>\in \aleph_{{\omega}+1}\times {\omega}_1$. 
We show that 
\begin{displaymath}\tag{$*_{{\alpha},j}$}
\mc H_{{\alpha},j}\subs \{a_{{\alpha},i}:i<j\}\cup
\bigcup\{\mc H_{{\xi},\ell}:{\xi}<{\alpha}\land a_{{\xi},\ell}\in F(a_{{\alpha},j})\}.
\end{displaymath}
Indeed, assume that 
$a_{{\eta},k}\in \mc H_{{\alpha},j}$. 
Observe that $
a_{{\eta},k}\subsetneq a_{{\alpha},j}$
implies 
$\<{\eta},k\>\triangleleft \<{\alpha},j\>$ by the definition of $\mc H$.

Since $\{a_{{\alpha},i}:i<j\}$ is a subset of the RHS of ($*_{}{\alpha},j$), we can assume that  ${\eta}<{\alpha}$. Then,  there is $z\in F(a_{{\eta},k})\cap F(a_{{\alpha},j})$
with $a_{{\eta},k}\subs z\subs a_{{\alpha},j}$ because $F$ is a wFN-function.
Hence $z=a_{{\xi},\ell}$ for some ${\xi}\le {\eta}$  and  $\ell<{\omega}_1$
by \ref{lm:wFN2enu}(\dag).
But $a_{{\eta},k}\in \mc H_{{\xi},\ell}$ and  $\mc H_{{\xi},\ell}$ is a subset of the 
RHS of $(**_{{\alpha},j})$.

So we proved $(*_{{\alpha},J  })$, which implies that $\mc H_{{\alpha},j}$ is 
countable by the inductive hypothesis. Thus Claim 2 holds. 

Hence $\mc H$ witnesses that (\ref{cs:wfn2nt}) holds.

To prove (\ref{cs:wfn2hnt}) we
recall a definition and a theorem  from \cite{ShSo23}.
   \begin{definition}[\cite{ShSo23}*{Def.~2.1}]
      We say that a family $\mc A\subs \br {\aleph_{\omega}};{\omega};$ is
      {\em nice }  if                                                                                                                                                                                                                                                                                                                               
      \begin{enumerate}[({N}1)]
      \item $\mc A$ is cofinal in ${[\aleph_{\omega}]}^{{\omega}}$, and 
      \item it is equipped with a well-ordering
      $\triangleleft$ of order type $\aleph_{{\omega}+1}$ such that 
      for each $B\in \mc A$ there is a countable subset 
      $$\mc I\subset\{A\in \mc A: A\triangleleft B\}$$
      with the following property: 
       for every $A\in \mc A$ with $A \triangleleft B$,
       there is  a finite subset $\mc J\subs \mc I$ such that  
      \begin{displaymath}
      A\cap B\subseteq^*    \bigcup \mc J.	
      \end{displaymath}
      \end{enumerate}
      \end{definition}
   Putting together Theorem 2.2 and Theorem 2.3 from \cite{ShSo23} we obtain:  
\begin{theorem}[\cite{ShSo23}]\label{tm:nice2group}
   If GCH holds and there is a nice family $\mc A\subs {[\aleph_{\omega}]}^{{\omega}}$, then $HnT$ holds.
\end{theorem}

   \smallskip
   \noindent{\bf Claim 3.} $\mc H$ is nice.

  Since  $\{a_{{\alpha},j}:{\alpha}<\aleph_{{\omega}+1},j<{\omega}_1\}$ is an
   enumeration of  ${[\aleph_{\omega}]}^{{\omega}}$ without repetition, 
   the lexicographical order $\triangleleft$ of $\aleph_{{\omega}+1}\times {\omega}_1$  induces a well-ordering  on 
   ${[\aleph_{\omega}]}^{{\omega}}$ of order type $\aleph_{{\omega}+1}$, which we also  denote by $\triangleleft$.

We show that the  $\triangleleft\restriction \mc H$ witnesses  (N2) for $\mc H$.

Let  $a_{{\alpha},j}\in \mc H$ be arbitrary.
Define
\begin{displaymath}
\mc F=\{z\in F(a_{{\alpha},j})\cap [a_{{\alpha},j}]^{\omega}:\exists {\gamma}_z<{\alpha}
\ \exists \ell_z<{\omega}_1\ z\subsq a_{{\gamma}_z,\ell_z}\in \mc H 
\},
\end{displaymath} 
and set 
\begin{displaymath}
\mc I=\{a_{{\gamma}_z,\ell_z}: z\in \mc F\}
\cup\big(\{a_{{\alpha},i}:i<j\}\cap \mc H\big). 
\end{displaymath}
We claim that $\mc I$ witnesses (N2) for $a_{{\alpha},j}$.
Clearly, $\mc I\subs \{a\in \mc H: a\triangleleft a_{{\alpha},j}\}$ is countable.

Assume that  $a_{{\eta},k}\in \mc H$ and   $a_{{\eta},k}\triangleleft a_{{\alpha},j}$, i.e.,     $\<{\eta},k\>\triangleleft \<{\alpha},j\>$.
If ${\eta}={\alpha}$, then $k<j$  and so
$$\mc J=\{a_{{\alpha},k}\}\in {[\mc I]}^{1} $$
satisfies the requirement of (N2): $a_{{\alpha},k}\cap a_{{\alpha},j}\subseteq a_{{\alpha},k}=\bigcup \mc J$.

Consider now the case  ${\eta}<{\alpha}$.
We may assume that  
   $x=a_{{\eta},k}\cap a_{{\alpha},j}$ is infinite. 
Since $F$ is a wFN-function and $x\subsq a_{{\alpha},j}$,    
    there exists $z\in F(x)\cap F(a_{{\alpha},j})$
   such that  $x\subseteq z\subsq a_{{\alpha},j}$.
Moreover, since $x\subseteq a_{{\eta},k}\subseteq A_{\eta}$ and $A_{\eta}$ is $F$-closed, it follows that 
$z\subs A_{\eta}$.  
Define 
\begin{displaymath}
\<{\gamma}^*,\ell^*\>=\min_{\triangleleft}\{\<{\gamma},\ell\>:z\subsq a_{{\gamma},\ell}\}.
\end{displaymath}
Then, $\gamma^*\le {\eta}$ and   $a_{{\gamma}^*,\ell^*}\in \mc H$, and so $\<{\gamma}^*,\ell^*\>$ witnesses that $z\in \mc F$.
Therefore, 
\begin{displaymath}
\text{$a_{{\gamma}_z,\ell_z}\in \mc I$ and 
   $a_{{\alpha},j}\cap a_{{\eta},k}=x\subsq z
\subsq a_{{\alpha},j}\cap a_{{\gamma}_z,\ell_z}$.}
\end{displaymath}
Thus, $$\mc J=\{a_{{\gamma}_z,\ell_z}\}\in {[\mc I]}^{1}$$
satisfies the requirements of (N2),  completing the proof of Claim 3. 

Hence, the family $\mc H$ is nice. So, by Theorem~\ref{tm:nice2group},
there exists an ${\omega}$-homogeneous, but not ${\omega}$-transitive
permutation group on $\aleph_{\omega}$. 
      \end{proof}

 \begin{problem} Are  there any additional implication  between the properties discussed in
   Section~\ref{sc:type}?
In particular, is it possible to reverse any of the implications shown  in Figure~\ref{fig:1}? 
\end{problem}     

Let us write {\bf (CCLC)} for the statement that  there exists a countably compact, locally countable regular space 
of size $\aleph_{\omega+1}$.  
Clearly, (SPL) implies (CCLC). By \cite{JuShSo88} and \cite{HaJuSo87}, if one iteratively
adds ${\omega}_1$ dominating reals  to a ground model, then in the resulting generic extension, 
both (SAT) and (CCLC) hold. 

This shows that  neither (SAT) not (CCLC) implies ($\neg$CC), because the partition relation 
$(\aleph_{{\omega}+1},\aleph_{{\omega}})\doublearrow  (\aleph_1,\aleph_0)$
is preserved by any c.c.c.\ forcing.

\begin{problem}
Is there any further implication between (CCLC) and the properties discussed in 
Section~\ref{sc:type}?
 In particular, is there any relation between 
(SAT) and (CCLC)?
\end{problem}

\bmhead{Acknowledgements}
The research on and preparation of this paper were supported by NKFI grant K129211.

\section*{Declarations}

\bmhead{Conflict of interest}
There is no competing interest to declare.

\end{document}